\documentclass[10pt,a4paper]{article}

\usepackage{amsfonts}
\usepackage{amsmath}
\usepackage{theorem}
\usepackage{amssymb}
\usepackage{epsfig}
\usepackage{subfigure}
\usepackage{enumerate}
\usepackage{color}
\usepackage[latin1]{inputenc}

\newcommand{\CC}{\mathbb{C}}

\newcommand{\RR}{\mathbb{R}}

\newcommand{\ints}{\int\limits}

\newcommand{\pa}{\partial}

\newcommand{\Om}{\Omega}

\title{A spectral approach to a constrained optimization problem for the Helmholtz equation in unbounded domains}

\author{Giulio Ciraolo\thanks{Corresponding author} \thanks{Dipartimento di Matematica e
Informatica, Universit\`a di Palermo, Via Archirafi 34, 90123, Italy.
\emph{E-mail}: {\tt g.ciraolo@math.unipa.it}}
\and Francesco Gargano\thanks{Dipartimento di Matematica e
Informatica, Universit\`a di Palermo, Via Archirafi 34, 90123, Italy.
\emph{E-mail}: {\tt francesco.gargano@unipa.it}}
\and Vincenzo Sciacca\thanks{Dipartimento di Matematica e
Informatica, Universit\`a di Palermo, Via Archirafi 34, 90123, Italy.
\emph{E-mail}: {\tt sciacca@math.unipa.it.}} }

\begin{document}

\maketitle

\begin{abstract}
We study some convergence issues for a recent approach to the problem of transparent boundary conditions for the Helmholtz equation in unbounded domains \cite{CGS13} where the index of refraction is not required to be constant at infinity. The approach is based on the minimization of an integral functional which arises from an integral formulation of the radiation condition at infinity.

In this paper, we implement a Fourier-Chebyshev collocation method to study some convergence properties of the numerical algorithm; 
in particular, we give numerical evidence of some convergence estimates available in literature \cite{Ci3}
and study numerically the minimization problem at low and mid-high frequencies. 
\end{abstract}

%
%
%
%
%
%



\section{Introduction}\label{intro}
We are concerned with time-harmonic wave propagation due to a source and some related classical scattering problem. The domain under consideration is $\RR^d,\ d=2,3$, wherein the field $u:\RR^d \to \CC$ is the \emph{outgoing
solution} of the Helmholtz equation (or \emph{reduced wave
equation})
\begin{equation}\label{Helmholtz intro}
\Delta u + k^2 n(x)^2u = f, \quad x \in \RR^d,
\end{equation}
where $k>0$ is the wavenumber, $n>0$ is the index of refraction, and $f$ is the source term. When $n$ is constant outside some compact region, say $n(x) \equiv 1$ for $x$ outside some large ball, the term \emph{outgoing solution}
means that we look for a solution $u$ of \eqref{Helmholtz intro} satisfying the following Sommerfeld radiation condition
\begin{equation} \label{rad cond Somm intro}
\lim_{r\to \infty} r^{\frac{d-1}{2}}\Big(\frac{\pa u}{\pa r} - iku \Big) =0,
\end{equation}
uniformly, where $r=|x|$. Condition \eqref{rad cond Somm intro} says that the radiation goes toward infinity and behaves at infinity as a spherical wave, decaying  uniformly in every direction.

For more general indexes of refraction $n$, \eqref{rad cond Somm intro} may be inadequate to guarantee the uniqueness of the solution for \eqref{Helmholtz intro}. As an example, one can think to an index of refraction which has some angular dependence, like $n(x)=\tilde{n}(x/|x|)$. In this case, it is known that the energy concentrates on lines rather than radiating in all directions (see \cite{PV} and references therein), and a pointwise condition like \eqref{rad cond Somm intro} seems to be not appropriate.

An extension of \eqref{rad cond Somm intro} to more general settings is given by the following condition
\begin{equation}\label{rad con Perth Vega intro}
\ints_{\RR^d} \Big| \nabla u(x) - ikn(x) u(x) \frac{x}{|x|} \Big|^2 \frac{dx}{1+|x|} < + \infty,
\end{equation}
see \cite{PV}. In \cite{PV}, it is proven that \eqref{rad con Perth Vega intro} guarantees the uniqueness of the outgoing solution of \eqref{Helmholtz intro} under very general assumptions on $n$. As an example, we mention that \eqref{rad con Perth Vega intro} applies to an index of refraction which satisfies $n(x)^2=n_\infty(x/|x|)^2 + p(x)$, where $n_\infty$ is a smooth angularly dependent function and $p$ is a perturbation with suitable decay at infinity (like $|x p(x)|\leq C$ for some constant $C>0$). For a correct and complete statement of the assumptions on $n$, we refer to \cite{PV}.
We remark that there are some extension of \eqref{rad con Perth Vega intro} to the Helmholtz equation with a magnetic potential (see \cite{Zu}). As far as the authors know, there are not results in this direction for the full Maxwell's equations.

A challenging issue in computational studies related to \eqref{Helmholtz intro} is how to deal numerically with the unbounded domain. Usually, one has to introduce a (bounded) computational domain $\Om$ and then prescribe boundary conditions on $\pa \Om$ which approximate the problem in the whole space. As it is well-known, a large amount of work has been done on this kind of problems. The most used methods are based on \emph{local or nonlocal conditions} involving $u$ (\cite{BT},\cite{EM},\cite{Gi4}), approximations of the \emph{Dirichlet-to-Neumann map} (\cite{KG},\cite{GiK},\cite{GrK}), \emph{infinite element schemes} (\cite{Be1}), \emph{boundary element methods} (see \cite{SS}) and the \emph{perfectly matched layer method} (\cite{Ber},\cite{TY}). For a deeper understanding of these problems and more recent developments, the interested reader can refer to the references cited in \cite{CGS13}.

Starting from \eqref{rad con Perth Vega intro}, in \cite{CGS13} we proposed a new approach to the computational study of the Helmholtz equation in unbounded domains. The idea is the following: we fix a computational bounded domain $\Om$ and approximate the solution of \eqref{Helmholtz intro} and \eqref{rad con Perth Vega intro} by the minimizer $u_\Om$ of the following constrained optimization problem:
\begin{equation}\label{introd minimiz problem PV}
\inf \{ J_\Omega[w]:\ w \textmd{ satifies } \Delta w + k^2 n(x)^2 w = f \textmd{ in } \Omega\},
\end{equation}
where
\begin{equation}\label{introd Funzionale}
J_\Omega[w]= \ints_{\Omega} \Big| \nabla w(x) - ikn(x) w(x) \frac{x}{|x|} \Big|^2  \frac{dx}{1+|x|}.
\end{equation}
In \cite{CGS13} we proved the existence and uniqueness of the minimizer $u_\Om$. Furthermore, $u_\Om$ converges to the outgoing solution in $H^1_{loc}$ norm as the computational domain $\Om$ enlarges and tends to cover the whole space. More precisely, if we consider $\Om=B_R$ (the ball or radius $R$ centered at the origin), then for any fixed $\rho >0$ we have that
\begin{equation*}
\| u_{B_R} - u\|_{H^1(B_\rho)} \to 0
\end{equation*}
as $R\to +\infty$; here, $u$ is the solution of \eqref{Helmholtz intro} satisfying \eqref{rad con Perth Vega intro}.

This approach is remarkable because: (i) it applies to a wide range of problems and for a large variety of indexes of refraction and source functions; (ii) numerical evidence (see \cite{CGS13}) shows that the geometry of the computational domain does not influence the accuracy of the numerical solution relevantly; (iii) it is suitable to be extended to other scattering problems (we have in mind the waveguide problem by using the results in \cite{AC},\cite{Ci1}--\cite{CM2}). On the other hand, the method has a slow convergence rate (in \cite{Ci3}, it is proved that $\| u_{B_R} - u\|_{H^1(B_\rho)} = O(R^{-1})$ as $R\to \infty$), and the computational cost is expensive (which is due to the volume integral formulation of the radiation condition). However, as already mentioned, this approach is of easy implementation in any setting and it is applicable even when the methods available in literature fail or are of difficult application. Indeed, most of those methods are based on the knowledge of the exact solution in the exterior of the computational domain or on a pointwise formulation of the radiation condition at infinity, which are not available -- for instance -- when $n$ is not constant at infinity. About the PML (perfectly matched layer) method, it is known that it can not be applied (at least in a standard way) when the index of refraction can not be continued analytically in the direction orthogonal to the artificial boundary (see \cite{OJ} and \cite{OZAJ}).

In this paper we shall investigate numerically how the numerical scheme depends on $R$ and $k$ and we will show some numerical examples which are of interest for the applications.

A relevant issue in this context is the study of the numerical algorithm at medium and high frequencies. When $k$ is large, the solution has a large amount of oscillations and the computational complexity of the problem increases. Indeed, it is well known that, for frequencies in the mid and high regime, the number of unknowns in the finite element methods scales at least like the cube of $k t$, where $t$ is a typical dimension of the scatterer \cite{IB}. In boundary integral methods the number of unknowns scales at least like the square of $kt$; this growth can be sometimes reduced by using some a priori knowledge regarding the oscillatory nature of the solution (see \cite{Gi} and reference therein).

From an analytic point of view, an argument as the one in \cite{Ci3} suggests that the convergence of the solution of \eqref{introd minimiz problem PV}-\eqref{introd Funzionale} to the exact solution improves for $k$ large. Since it is difficult to obtain an exact estimate of the rate of convergence, one of the goals of this paper is to test numerically the convergence of the algorithm for $k$ large. The numerical results presented in this paper (for $n(x)\equiv 1$) suggest that our approach gives consistent results at mid-high frequencies. However, since our approach has a large computational cost, we can not perform numerical simulations for $k$ very large (in our tests, $k$ varies up to 96).

In this paper we use a spectral collocation method to implement the numerical simulations and calculate the minimizer of \eqref{introd minimiz problem PV}. Spectral methods are largely employed in problems where the solution is known to be analytic, and they usually provide a fast convergence of the numerical scheme. However, the solution of a typical scattering problem is expected to be no more than $C^{1,\alpha}$ regular (the coefficient $n$ may be just a $L^\infty$ function). The main reason for using a spectral method here is due to the \emph{global} nature of the method. Indeed, since we do not have to prescribe boundary conditions to the equation but we have to minimize an integral functional which is defined on the whole computational domain, the implementation by using a global method may reduce the computational complexity of the numerical scheme.
We believe that, since the numerical studies presented in this paper will not make use of any approximation which could be used for $k$ large (see for instance \cite{Gi}), the results presented in this paper are relevant in the case of mid-high frequency regime.


The spectral approach will be used also to implement numerical simulations at low frequencies ($k\approx 0$). In this case, it is expected that the error norms are worse as $k$ becomes smaller. This is confirmed by some numerical simulations presented in this paper. As we will show numerically, it is interesting to notice that this behavior does not appear for any choice of the source function.

The paper is organized as follows. In Section \ref{sect 2} we recall some analytical results and we conjecture the behavior of the numerical scheme for frequencies in the low, mid and high frequency regime. The study is presented in two dimensions but it can be extended to more dimensions. In Section \ref{numeric_scheme} we describe the numerical scheme and the spectral collocation method that we used. In Section \ref{numerical_results} we study the source problem \eqref{introd minimiz problem PV}-\eqref{introd Funzionale} and present some numerical results that confirm our conjectures on the convergence rates. Moreover, some numerical simulations for an angular dependent index of refraction are also shown. In Section \ref{sect 5} we show numerical simulations for some classical scattering problem.


\section{Preliminaries and general considerations} \label{sect 2}
Let $u$ be the solution of \eqref{Helmholtz intro} satisfying \eqref{rad con Perth Vega intro} and denote by $B_R$ the ball of radius $R$ centered at the origin. Let $u_R$ be the solution of the following minimization problem
\begin{equation}\label{minimiz problem PV sect 2}
\inf \{ J_{B_R}[w]:\ w \textmd{ satifies } \Delta w + k^2 n(x)^2 w = f \textmd{ in } B_R\},
\end{equation}
where $J_{B_R}$ is given by \eqref{introd Funzionale} with $\Om=B_R$. In \cite{CGS13} and \cite{Ci3}, it is proved that $u_R$ approximate $u$ in $H^1_{loc}$ norm, more precisely that, for any fixed $\rho>0$, we have that
\begin{equation}\label{sec2 conv rate}
 \|u_R - u\|_{H^1(B_\rho)} \approx \frac{1}{R}, \quad \textmd{as } R \to +\infty,
\end{equation}
for $k$ fixed. As it was noticed in \cite{CGS13}, there is an interesting interplay between the parameters $k,\: R $ and $f$ in \eqref{sec2 conv rate}. In this paper, we shall investigate  \eqref{sec2 conv rate} and some aspects of its dependenc on these parameters. In particular, by using numerical results in the simplest case possible ($n(x) \equiv 1$ and $d=2$), we are interested to study how the convergence in \eqref{sec2 conv rate} depends on the wavenumber $k$.

In \cite{CGS13}, it was noticed that the convergence is expected to worsen as $k$ becomes smaller (See Remark 3.1 in \cite{CGS13}), that is
\begin{equation} \label{rate conv k0}
\lim_{k\to 0^+} \|u_R - u\|_{H^1(B_\rho)} = +\infty, \quad R \textmd{ fixed}.
\end{equation}
It is not surprising that there is a behaviour like \eqref{rate conv k0} for $k$ small. Indeed, when $k=0$ the Helmholtz equation reduces to the Laplace equation and \eqref{rad con Perth Vega intro} (or the Sommerfeld radiation condition \eqref{rad cond Somm intro}) is not appropriate to guarantee that the scattering problem is well-posed. The first two goals of this paper are: (i) a deeper study of the rate of blow-up in \eqref{rate conv k0} and (ii) to understand if there are some sources for which there is no blow-up of $\|u_R - u\|_{H^1(B_\rho)}$ as $k\to 0^+$.

Now, let us discuss the behavior of the minimization problem \eqref{minimiz problem PV sect 2} in the high-frequency regime. Following Remark 3.1 in \cite{CGS13}, it is conjectured that the convergence improves for $k$ large:
\begin{equation} \label{rate conv kinfty}
\lim_{k \to \infty} \|u_R - u\|_{H^1(B_\rho)}=0, \quad  R \textmd{ fixed}.
\end{equation}
We do not know the exact rate of the convergence in \eqref{rate conv kinfty}: a deeper study of this issue is another goal of this paper.

To achieve this goal, it is needed to study the Helmholtz equation numerically for large values of $k$. Due to the large number of oscillations, this is a nontrivial issue: it is known that the computational complexity needed to solve the Helmholtz equation with some boundary condition (and by using a finite element scheme) grows as $k^3$ for $k\to \infty$ (see \cite{Gi} and references therein). Here, the implementation of \eqref{minimiz problem PV sect 2} requires a larger computational effort. Indeed, since we have to minimize a volumetric integral functional subject to a linear constrain, by Lagrange multipliers we have to deal with the inversion of a linear system which has a larger dimension and the computational complexity will increase at least as $k^6$.

For this reason, it is reasonable to use spectral methods to implement the
numerical algorithm. Indeed, a spectral method is a global method and it is expected
that less points are needed to obtain the desired accuracy of the numerical algorithm, specially when
the index of refraction and the source function are analytical\footnote{We implemented \eqref{minimiz problem PV sect 2}
by using a finite element scheme and we noticed that the grid points needed in this case were $10^2$ times more than in the spectral collocation
method ($10^6$ and $10^4$ grid points, respectively).}. We stress that the spectral method is still applicable when $n$ and $f$ are not
continuous functions, which implies that the second derivatives
of the solution may be discontinuous. However the solution will be at least of class $C^1$ and it is well known that
spectral methods can be applied to solve the numerical problem (see for instance \cite{quart}).


\section{Numerical scheme}\label{numeric_scheme}

We use a spectral Fourier-Chebyshev collocation method to implement the constrained optimization problem
\eqref{introd minimiz problem PV}.
\footnote{
We notice that, thanks to the collocation method,
the accuracy of a Fourier-Chebyshev-type scheme is equivalent to
the one of a Fourier-Jacobi method when one uses (roughly) the same order of polynomials \cite{LJW,MC,VI,VII,MM,Boy01}.

Collocation methods for ODEs and PDEs can also be implemented by using different basis, like Hermite or Legendre polynomials (\cite{OG,KS}); however, those approaches require a suitable decay at infinity which is not satisfied in our case (see also \cite{SW} for a more exhaustive reading
on the recent advances on the spectral methods).
}

We implement the minimization problem for $\Om = B_R$ in $\RR^2$, a disk of radius $R$
centered at the origin which is parameterized in terms
of the usual polar coordinates $(\rho,\theta)$.
To guarantee the well-posedness of the numerical scheme at $\rho=0$, we follow \cite{For95,For96,Tref01} and consider the symmetry conditions $u
(\rho,\theta) = u(-\rho,(\theta + \pi) \mod 2\pi )$ and $n
(\rho,\theta) = n(-\rho,(\theta + \pi) \mod 2\pi )$ in the
$(\rho,\theta)$-space, which implies that we look at the solution in the \textit{extended}  domain $\Om^E=[-R,R]\times(0,2\pi]$.

We discretize $\Om^E$ by using a periodic Fourier grid in the angular variable $\theta$ and Fourier
points $\theta_i=\frac{ 2\pi i}{M_\theta}$, $i=1,\ldots, M_\theta$ ($M_{\theta}$ taken odd). We use a Chebyshev grid in the radial variable $\rho$
with Gauss-Lobatto points $\rho_j=R \cos(\frac{j\cdot \pi}{2M_\rho-1})$, $j=0, \ldots, 2M_\rho-1$, so that $\rho=0$ is not included in the computational grid.

We write the constraint \eqref{introd minimiz problem PV} in $\Om^E$ in polar coordinates
$$ u_{\rho \rho} + \frac{1}{\rho} u_{\rho} +\frac{1}{\rho^2} u_{\theta \theta} + k^2 n(\rho, \theta)^2 u
= f\label{constrain_polar}
$$
and we discretize the Laplacian by the well-known differential matrices for
spectral collocation methods (see \cite{Boy01,For96,Tref01}). This means that, at any collocation
grid points, we can construct differentiation formulas for the derivatives of $u$ in term of the values of $u$ itself at all collocation points. We define the
matrices
$\mathcal{D}^{rr}=(d^{\rho\rho})_{l,m},\mathcal{D}^{r}=(d^{\rho})_{l,m},\mathcal{D}^{\theta}=(d^{\theta})_{p,q},\mathcal{D}^{\theta\theta}=(d^{\theta\theta})_{p,q},\mathcal{T}=(t)_{w,z}$, where
\begin{eqnarray*}
&& d^{\rho}_{l,m}=\frac{\bar{c}_l(-1)^{l+m}}{\bar{c}_m(\rho_l-\rho_m)}, \quad 0\leq l,m,\leq 2M_{\rho}-1, l\neq m\\
&& d^{\rho}_{l,l}=-\sum_{m\neq j} d^{\rho}_{l,m},\quad 0\leq l \leq 2M_{\rho}-1,\\
&& d^{\rho\rho}_{l,m}=\sum_{k=0}^{2M_{\rho}-1}d^{\rho}_{l,k}d^{\rho}_{k,m},
\end{eqnarray*}
for the radial derivatives,
\begin{eqnarray*}
&& d^{\theta}_{p,q}=\frac{(-1)^{p+q}}{2\sin((\theta_{p+1}-\theta_{q+1})/2)},\quad 0\leq p,q\leq M_{\theta}-1, p\neq q,\\
&& d^{\theta}_{p,p}=0,\quad0\leq p\leq M_{\theta}-1,\\
&& d^{\theta\theta}_{p,q}=(-1)^{p+q+1}\frac{\cos((\theta_{p+1}-\theta_{q+1})/2)}{2\sin^2((\theta_{p+1}-\theta_{q+1})/2)},\quad 0\leq p,q\leq M_{\theta}-1, p\neq q,\\
&& d^{\theta\theta}_{p,p}=-\frac{(M_{\theta}+1)^2-1}{12},\quad0\leq p\leq M_{\theta}-1,
\end{eqnarray*}
for the angular derivatives, and
\begin{eqnarray*}
&& t_{w,w}= n(\rho_{w -M_{\rho}\cdot[w/M_{\rho}]},\theta_{[w/M_{\rho}]+1}),\quad 0\leq w\leq M_{\rho}\cdot M_{\theta}-1,\\
&& t_{w,z}=0,\quad 0\leq w,z\leq M_{\rho}\cdot M_{\theta}-1,w\neq z,
\end{eqnarray*}
where $\bar{c}_0=\bar{c}_{(2M_{\rho}-1)}=2,\bar{c}_j=1$ for $j\neq 0,(2M_{\rho}-1)$, and $[\cdot]$ is the integer--part function.
The matrices $D^{r}$ and $D^{rr}$ have the general form
$$
\mathcal{D}^{\sharp}= \left( \begin{array}{cc}
         \mathcal{D}^{\sharp}_1 & \mathcal{D}^{\sharp}_2\\
        \mathcal{D}^{\sharp}_3 & \mathcal{D}^{\sharp}_4
        \end{array}\right)
$$
where $D^{\sharp}_1,D^{\sharp}_2$ and $D^{\sharp}_3,D^{\sharp}_4$ are the $M_{\rho}\cdot M_{\rho}$ blocks relative to the subregions of $\Om^C$ in which $\rho>0$ and $\rho<0$, respectively.
Due to the symmetry condition
$$u(\rho,\theta) = u(-\rho,(\theta + \pi) \mod 2\pi ),$$
the contributions due to the lower blocks of the matrices $\mathcal{D}^{\sharp}$ are redundant,
and therefore we neglect  these contributions.
Hence, we can consider only the values corresponding to $\rho>0$, and we define the following complex-valued vectors
\begin{eqnarray*}
v_{j+M_\rho (i-1)}&=&u(\rho_j,\theta_i),\\
\phi_{h+M_\rho (h-1)}&=&f(\rho_h,\theta_i),
\end{eqnarray*}
for $i=1,..,M_\theta$, $j=0,..,M_\rho-1$ and $h=1,..,M_\rho-1$.

We stress that the minimizer of Problem \ref{introd minimiz problem PV}
has to satisfy the Helmholtz equation only in the interior points of $\Om^E$. For this reason, the $M_\theta$ boundary values of $u$ in $\rho=R$
are the degrees of freedom of the problem and we do not need to evaluate $f$ at the boundary points.    Therefore the
discretized Helmholtz equation can be seen as an algebraic constrain of the form
\begin{equation}\label{Helm_discreto}
\tilde A v = \phi,
\end{equation}
where $\tilde A$ is the $M_\theta\cdot (M_\rho-1) \times M_\theta\cdot M_\rho$ matrix obtained by removing the first row of the matrix
\begin{eqnarray*}
A=(\mathcal{D}^{rr}_1+\mathcal{R}  \mathcal{D}^{r}_1)\otimes \left(\begin{array}{cc} I_{(M_{\theta}-1)/2} & 0\\ 0 & I_{(M_{\theta}+1)/2}\end{array} \right)+\\
(\mathcal{D}^{rr}_2+\mathcal{R}  \mathcal{D}^{r}_2)\otimes \left(\begin{array}{cc} 0 & I_{(M_{\theta}+1)/2}\\ I_{(M_{\theta}-1)/2} & 0\end{array} \right)+\mathcal{R}^2
\otimes \mathcal{D}^{\theta\theta}+k^2\mathcal{T},
\end{eqnarray*}
being $\mathcal{R}=\text{diag}(\rho_j^{-1})$, $0\leq j\leq M_{\rho}-1$
and $I_j$ the identity matrix of order $j$.

The functional \eqref{introd Funzionale} is written in polar coordinates as \begin{eqnarray*}\label{Funzionale_polar}
J_\Omega[u]= \ints_{0}^{R} \ints_{0}^{2\pi} \left( \left| \cos(\theta)
u_\rho - \frac{\sin(\theta)}{\rho} u_\theta  -ik \cos(\theta) n(\rho, \theta)
u\right|^2 + \right. \\ + \left. \left| \sin(\theta) u_\rho +
\frac{\cos(\theta)}{\rho} u_\theta  -ik \sin(\theta) n(\rho, \theta) u
\right|^2 \right) \frac{\rho}{1+\rho} d \rho d \theta,
\end{eqnarray*}
and we use Gauss quadrature formula as
integration rule in the functional. Therefore the discretized functional can be written as a quadratic form:
\begin{equation}\label{Funzionale_discreto}
J_{B_R}[v]= \frac{1}{2} \bar{v}^T \tilde H v,
\end{equation}
where the matrix $\tilde H=H^{\dagger}\mathcal{W}H$ has dimension
$M_\theta\cdot M_\rho\times M_\theta\cdot M_\rho$, being

$$H=\left( \begin{array}{c}\mathcal{D}^{r}_1\otimes \left(\begin{array}{cc} I_{(M_{\theta}-1)/2} & 0\\ 0 & I_{(M_{\theta}+1)/2}\end{array} \right)+
\mathcal{D}^{r}_2\otimes \left(\begin{array}{cc} 0 & I_{(M_{\theta}+1)/2}\\ I_{(M_{\theta}-1)/2} & 0\end{array} \right)-ik\mathcal{T}\mathcal{C} \\
\\
\mathcal{R}\otimes D^{\theta} -ik\mathcal{T}\mathcal{S}\end{array}\right),$$
$\mathcal{C}=\text{diag}(\cos(\theta_{[j/ M_{\rho}]+1}))$, $0\leq j\leq  M_{\rho}\cdot M_{\theta}-1$, $\mathcal{S}=\text{diag}(\sin(\theta_{[j/ M_{\rho}]+1}))$, $0\leq j\leq  M_{\rho}\cdot M_{\theta} -1$,

$$\mathcal{W}=\left(\begin{array}{cc}\mathcal{W}_1 & 0\\ 0 & \mathcal{W}_2 \end{array} \right),$$
with
\begin{eqnarray*}
 \mathcal{W}_1=\mathcal{W}_2=\text{diag}(\frac{2R\pi\theta_1(1-\rho_{j-M_{\rho}\cdot[j/M_{\rho}]}^2)^{-1/2}}{\tilde c_j (2M_{\rho}-1)}\cdot\frac{\rho_{j-M_{\rho}\cdot[j/M_{\rho}]}}{1+\rho_{j-M_{\rho}\cdot[j/M_{\rho}]}} ),\\
  0\leq j\leq M_{\rho}\cdot M_{\theta}-1,\quad c_0=2, \quad c_j=1 \text{ for } j\geq0
\end{eqnarray*}
The matrix $H$  determines the components of the vector ($\nabla u-iknu\frac{x}{|x|}$) in polar coordinates,
while the matrix $\mathcal{W}$ is obtained by the Gauss quadrature formula derived from the Chebyshev approximatin in the Gauss-Lobatto points (see \cite{pey})
Hence, the discretized problem can be stated as follows: to find a complex-valued vector $v$ of size $M_\theta\cdot M_\rho$ which minimizes the problem
\begin{equation} \label{discretized pb}
\min \frac{1}{2} \bar{v}^T \tilde H v \textmd{ such that } \tilde A v= \phi.
\end{equation}

As it is well-known, constrained optimization problems involve a set of Lagrange
multipliers $\lambda$. By standard optimization theory, we have that the minimizer of \eqref{discretized pb} is given by the solution of
the following (sparse) algebraic linear system:

\begin{equation}\label{sist_lagrange}
\left(
\begin{array}{cc}
\tilde A & 0 \\
\tilde H & \tilde A^T
\end{array}
\right)
\left(
\begin{array}{c}
v \\
\lambda
\end{array}
\right)
=
\left(
\begin{array}{c}
\phi \\
0
\end{array}
\right),
\end{equation}
where the vector of Lagrange multipliers $\lambda$ has dimension $M_\theta\cdot (M_\rho-1)$.

We compute the solution of \eqref{sist_lagrange} by a least square iterative method (choosing a tolerance of
$10^{-10}$).
To test our numerical scheme, we assume that $n(x)=1$. In this case, the exact solution $u$ is given by
\begin{equation*}
u(x)=\frac{i}{4}\int_{\RR^2} H_0^{(1)}(k|x-y|) f(y) dy,
\end{equation*}
where $H_0^{(1)}(\cdot)$ is the zeroth-order Hankel function of the first kind.
As source term in \eqref{Helmholtz intro}, we consider
\begin{equation}\label{tipo_sorg_0}
f(x)= - e^{-\sigma |x|^2},
\end{equation}
with the choice of $\sigma=30$.

We tested the convergence properties of the numerical scheme by evaluating several error norms ($L^2$, $H^1$ and $L^\infty$ as in \cite{CGS13}) and observed the same order of errors as in \cite{CGS13}. However, thanks to the spectral approach, we can use less grid
points. We checked that the numerical solution does not depend on the spatial resolution that we choose (both in the angular and radial variables). For instance, by fixing $R=4$, $k=1$ and  $M_\rho=100$, the
errors remain at order
$10^{-3}$ for values of $M_\theta$ from $11$ to $81$.
Viceversa, analogous results can be observed if we fix $M_\theta=21$
and let $M_\rho$ vary from $50$ to $600$. It results that we obtain satisfying numerical results by using a not too high resolution and, as a consequence,
we can perform numerical simulations for large values of $R$ and $k$.

\vspace{0.5em}

We conclude this section by noticing that if one considers the 3-D case, the numerical approach is very similar. For instance, one can use a Fourier grid in each angular variable and a Chebyshev grid in the radial direction.

\section{Numerical Results I: source problems}\label{numerical_results}

In this section we study some convergence properties for the optimization problem \eqref{discretized pb}. In particular, we show numerical evidence for the convergence rate estimate \eqref{sec2 conv rate}, and we study how the numerical algorithm behaves at low and mid/high frequencies. Moreover, we present some numerical simulations for a problem where the index of refraction is an unbounded perturbation of an angular dependent background function.

\subsection{Convergence results}
We perform our analysis by choosing a constant index of refraction $n(x)=1$ and for the following choices of the source term $f$:
\begin{eqnarray}
f_1(x) & = &  e^{-\sigma|x|^2}, \label{tipo_sorg_1}\\
f_2(x) &=& e^{-\sigma|x-p_+|^2} - e^{-\sigma|x-p_-|^2}, \label{tipo_sorg_2}\\
f_3(x) & = & e^{-\sigma|x-p_+|^2} - e^{-\sigma|x-p_-|^2}+e^{-\sigma|x-q_+|^2}, \label{tipo_sorg_3}\\
f_4(x) &=& e^{-\sigma|x-p_+|^2} - e^{-\sigma|x-p_-|^2}+e^{-\sigma|x-q_+|^2} - e^{-\sigma|x-q_-|^2}, \label{tipo_sorg_4}
\end{eqnarray}
where $p_\pm = \left(\pm 0.25,0 \right), q_\pm = \left(0,\pm 0.25 \right)$  and with $\sigma = 30$. The
computational domain is the disk $B_R$, where $R$ will be specified from time
to time.

The choice of a simple index of refraction and of quite simple analytic source terms is motivated by the fact the here we want to study the behavior of the minimization problem for $R$ and $k$ large, and such choices of $n$ and $f$ permits to consider larger values of $k$ and $R$ (the spectral method does its very best with such choices).

\paragraph{Dependence on $R$.}
We give numerical evidence of the estimate \eqref{sec2 conv
rate}. We fix $\rho=1$ and compute the numerical solution for several values of $R$, which range from $R=4$  to $R=32$. For a fixed $R$, the solution is obtained
by taking $M_{\theta}=21$ and $M_{\rho}=25 R$.

To obtain the estimate \eqref{sec2 conv
rate}, we need to interpolate the solution $u_{R}$ of
problem \eqref{discretized pb} -- which is defined in the ball $B_R$ -- in
$B_1$. The interpolation is performed by using the spectral Fourier-Chebyshev expansion of the solution:
firstly, we compute the Fourier-Chebyshev modes
by writing the discrete Fourier-Chebyshev transform of the solution $u_R$ in $B_R$, that is
$$
u_R(\rho,\theta)=\sum_{k=-(M_{\theta}-1)/2}^{{(M_{\theta}+1)/2}} \sum_{h=0}^{2M_{\rho}-1} u_{hk}e^{-ik\theta}T_h(\rho),
$$
where $T_h(\rho)$ is the Chebyshev $h-$polinomial of the first kind.
Then, we evaluate the Fourier-Chebyshev expansion in the points of
$B_1$, which is discretized by using $M^{B_1}_{\theta}=21$ and $M^{B_1}_{\rho}=100$ grid points in angular and radial variables respectively.

\begin{figure}
\begin{center}
\vspace*{-1.8cm}\includegraphics[width=\columnwidth]{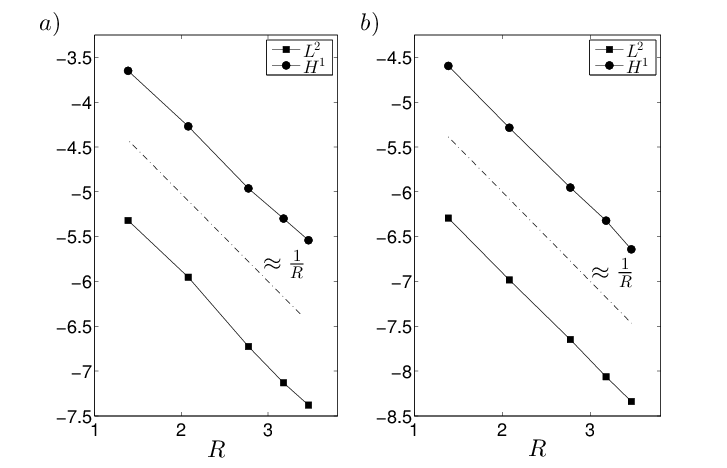}
\caption{The values of the error norms in $B_1$ for $f=f_1$ ($\bold{a})$ and for $f=f_2$ ($\bold{b})$ in the log-log scale. Here, $k=1$, $M^{B_1}_{\theta}=21$ and $M^{B_1}_{\rho}=100$.}
\label{stima_decrescita_tab_1}
\end{center}
\end{figure}

In Figures \ref{stima_decrescita_tab_1}a)-b)  we show the error norms in $B_1$ in the log-log scale for $f=f_1$ and $f=f_2$, respectively. We notice that the rate of decay is proportional to $1 / R$, which perfectly agrees with \eqref{sec2 conv rate}.
%
%

\paragraph{Dependence on $k$.}
In this subsection we study numerically how the optimization problem behaves for
small and large values of $k$. We fix the computational domain $B_R$ with $R=4$,
and we compute the numerical solution for several values of $k$, starting from
$k=1/32$ and up to $k=96$. The resolution is fixed by taking
$M_\theta=21$ and $M_\rho=800$.

\begin{figure}
\begin{center}
\vspace*{-1cm}\includegraphics[width=\columnwidth]{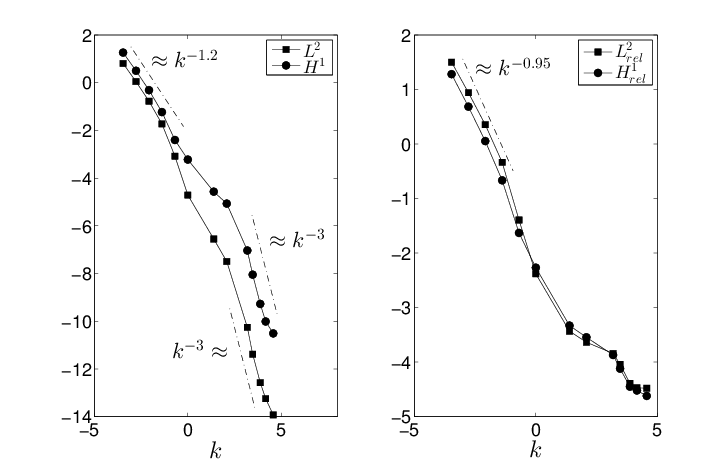}
\caption{The values of the error norms at several values of $k$ for $f=f_1$ in the log-log scale. Here, $R=4$,
$M_\theta=21$,
$M_\rho=800$.}
\label{fig_k_mono}
\end{center}
\end{figure}

\begin{figure}
\begin{center}
\vspace*{-2cm}\includegraphics[width=\columnwidth]{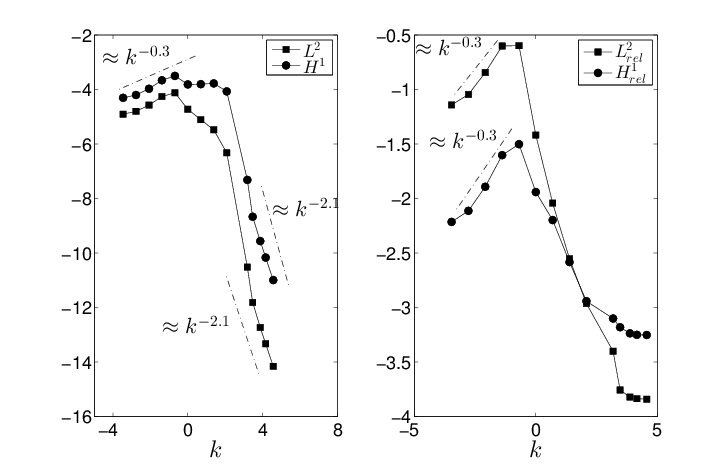}
\caption{The values of the error norms at several values of $k$ for $f=f_2$ in the log-log scale. Here, $R=4$,
$M_\theta=21$,
$M_\rho=800$.}
\label{fig_k_di}
\end{center}
\end{figure}

\begin{figure}
\begin{center}
\vspace*{-1.8cm}\includegraphics[width=\columnwidth]{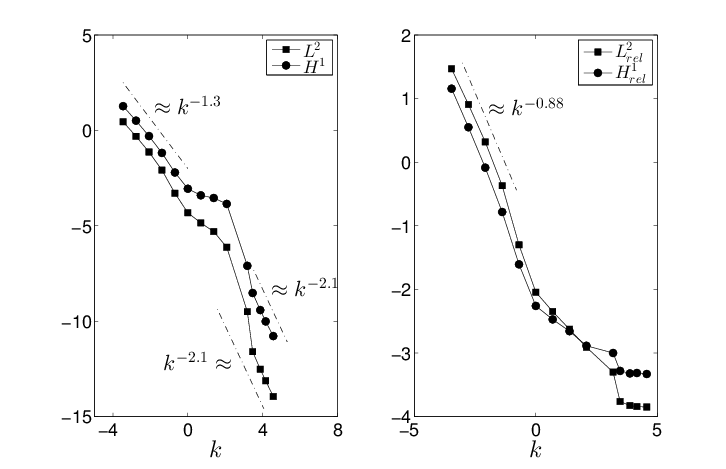}
\caption{The values of the error norms at several values of $k$ for $f=f_3$ in the log-log scale. Here, $R=4$,
$M_\theta=21$,
$M_\rho=800$.}
\label{fig_k_tri}
\end{center}
\end{figure}

\begin{figure}
\begin{center}
\vspace*{-2cm}\includegraphics[width=\columnwidth]{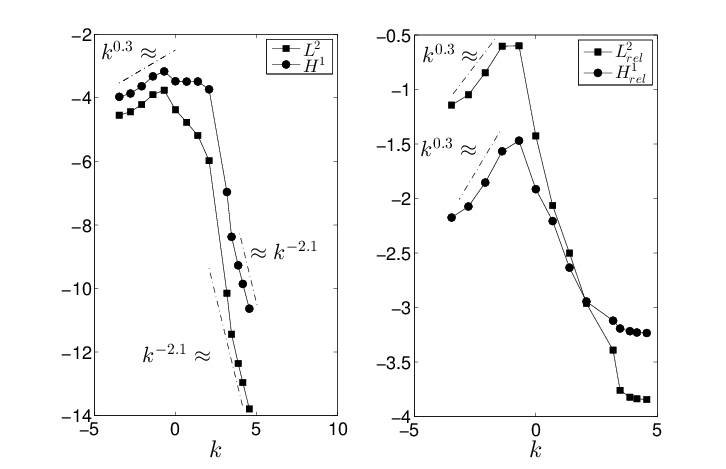}
\caption{The values of the error norms at several values of $k$ for $f=f_4$ in the log-log scale. Here, $R=4$,
$M_\theta=21$,
$M_\rho=800$.}
\label{fig_k_quadri}
\end{center}
\end{figure}

In Figs.\ref{fig_k_mono}-\ref{fig_k_quadri} we show the
values of the $L^{2}$,$L^2_{rel}$, $H^1$ and $H^1_{rel}$ error norms in the log-log scale for the sources term specified in \eqref{tipo_sorg_1}-\eqref{tipo_sorg_4}, respectively.
We observed a different behavior according to the choice of the source term. In cases $f=f_1$ and $f=f_3$  the error
norms deteriorates considerably as $k$ becomes smaller. For $f=f_2$ and $f=f_4$, the error norms still deteriorates as $k\to 0$ but with a slower rate.

When $k$ becomes larger, Figures \ref{fig_k_mono}-\ref{fig_k_quadri} show that the error norms improve in all cases.
The numerical simulations suggest that the $L^2$ and $H^1$ convergence improves
as $k^{-\alpha}$, where $\alpha \approx 3$ for $f_1$ and $\alpha \approx 2.1$ for $f_2,f_3$ and $f_4$. We also notice that the $L^2_{rel}$ and $H^1_{rel}$ error norms are stable for $k$ large. This behavior is remarkable because it give numerical evidence that the numerical algorithm is stable in the mid frequency regime. This is in accordance with analytic studies for $k$ large. However, due to computational limits, we did not test the method at (very) high frequencies. Since the computational grid must be chosen finer and finer as $k\to \infty$, the numerical results deteriorate for $k$ very large.


\subsection{Background medium with variable index of refraction.}
We assume that the background medium has variable index of refraction
$$n_0(x)^2= 2+\frac{x_1}{|x|},$$
and that the source has compact support. In particular, we consider the following two source problems:
\begin{description}
 \item[PVR 1]: $\Delta u + k^2 n(x)^2 u = \chi_{Q_{0.5}(0,0)}(x)$,
 \item[PVR 2]: $\Delta u + k^2 n(x)^2 u = \chi_{Q_{0.5}(0.5,0)}(x)-\chi_{Q_{0.5}(-0.5,0)}(x)$,
\end{description}
where
\begin{equation*}
n(x)^2= n_0(x)^2 + e^{-|(x_1-1,x_2)|^2},
\end{equation*}
and $Q_r(p)$ is the square centered at $p\in\RR^2$ of side $2r$.

\begin{center}
\begin{figure}
\hspace*{-1.2cm}\subfigure[Problem \textbf{PVR 1}, real part: $k=1$]{\includegraphics[width=7.3cm]{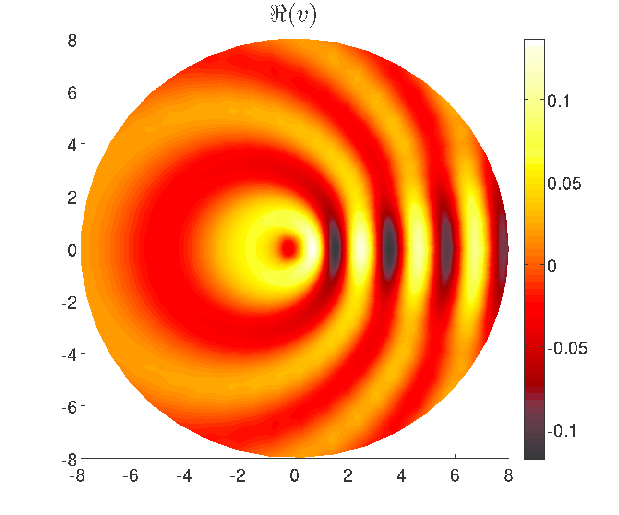}}
\hspace*{-0.5cm}\subfigure[Problem \textbf{PVR 1}, imaginary part: $k=1$]{\includegraphics[width=7.3cm]{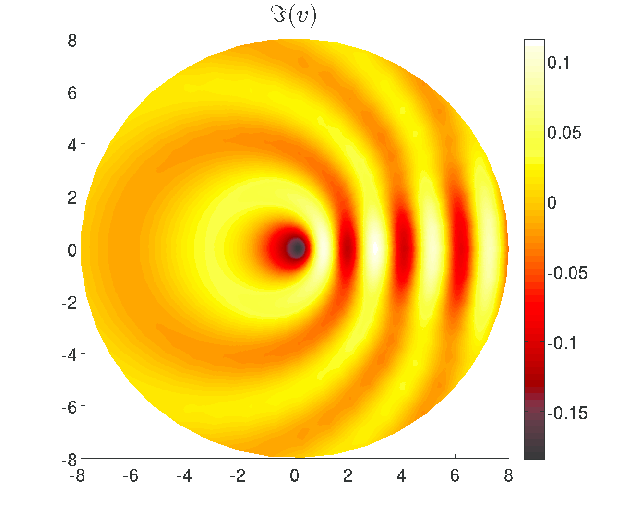}}
\hspace*{-1.2cm}\subfigure[Problem \textbf{PVR 1}, real part: $k=2$]{\includegraphics[width=7.3cm]{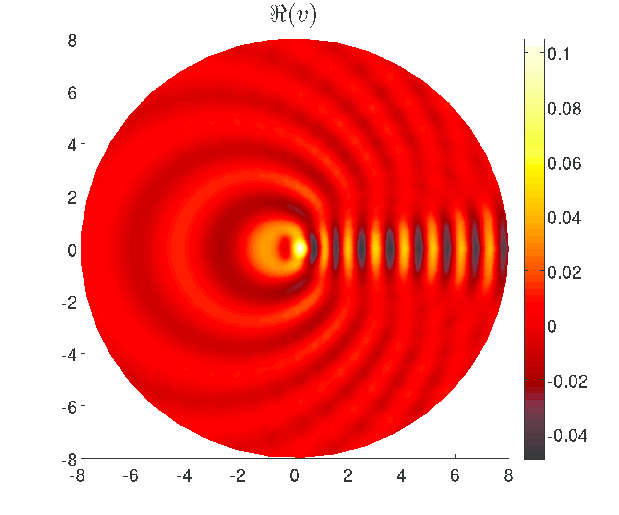}}
\hspace*{-0.5cm}\subfigure[Problem \textbf{PVR 1}, imaginary part: $k=2$]{\includegraphics[width=7.3cm]{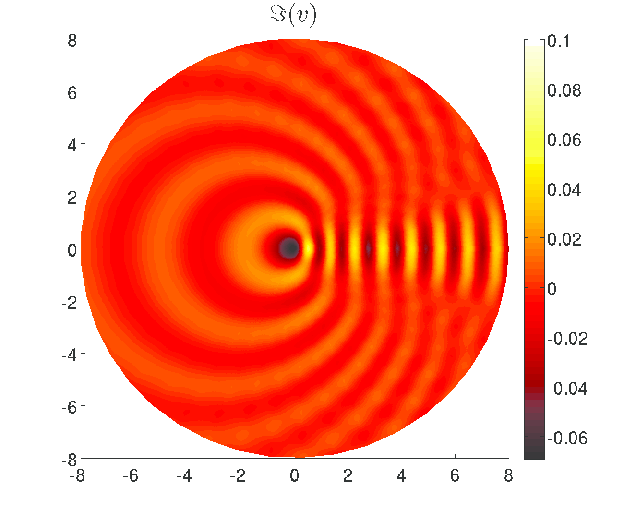}}
\caption{\textbf{Problem PVR 1:} real and imaginary part of the numerical solution for $k=1,2$.}
\label{ppvr1}
\end{figure}
\end{center}

\begin{center}
\begin{figure}
\hspace*{-1.2cm}\subfigure[Problem \textbf{PVR 2}, real part: $k=1$]{\includegraphics[width=7.3cm]{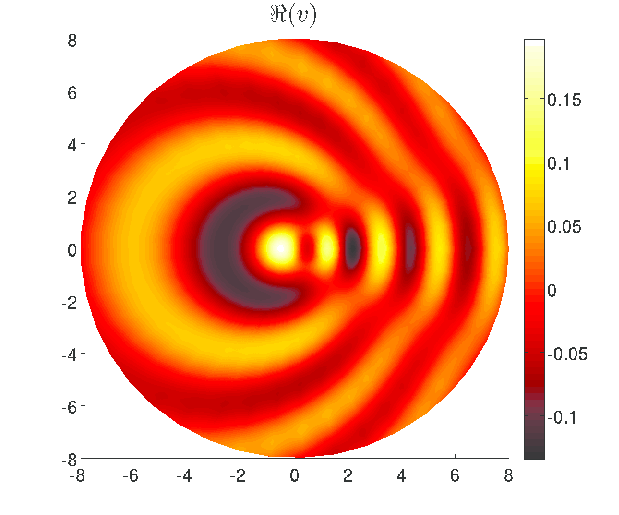}}
\hspace*{-0.5cm}\subfigure[Problem \textbf{PVR 2}, imaginary part: $k=1$]{\includegraphics[width=7.3cm]{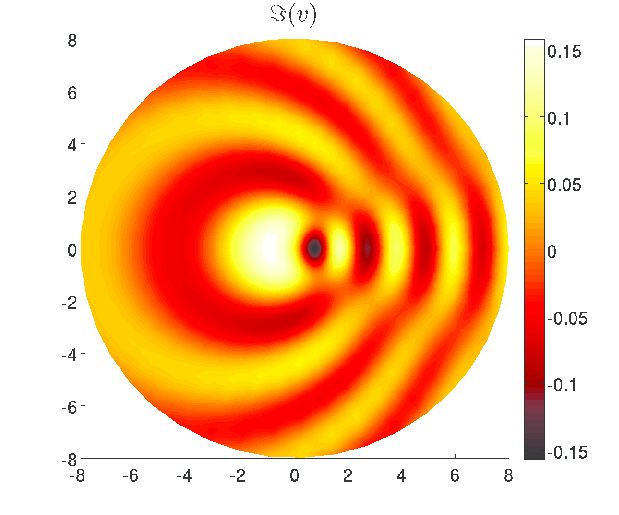}}
\hspace*{-1.2cm}\subfigure[Problem \textbf{PVR 2}, real part: $k=2$]{\includegraphics[width=7.3cm]{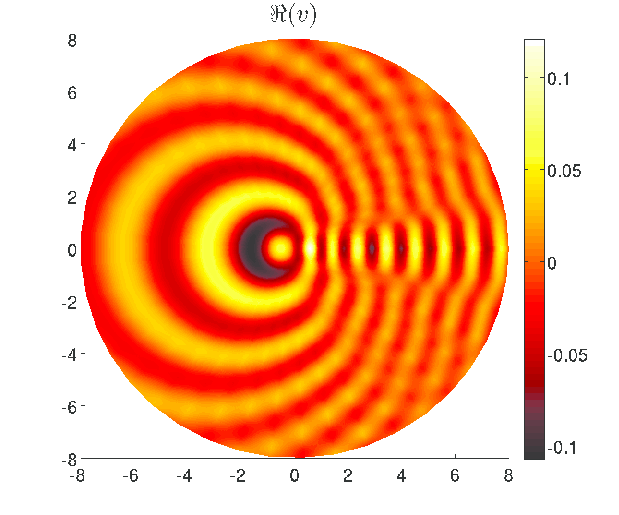}}
\hspace*{-0.5cm}\subfigure[Problem \textbf{PVR 2}, imaginary part: $k=2$]{\includegraphics[width=7.3cm]{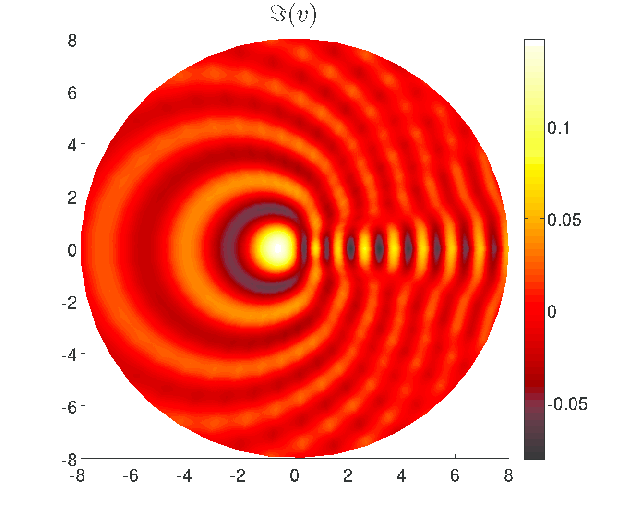}}
\caption{\textbf{Problem PVR 2:} real and imaginary part of the numerical solution for $k=1,2$.}
\label{ppvr2}
\end{figure}
\end{center}

\begin{center}
\begin{figure}
\hspace*{-1.2cm}\subfigure[Problem \textbf{PVR 1}: $k=1$]{\includegraphics[width=7.5cm]{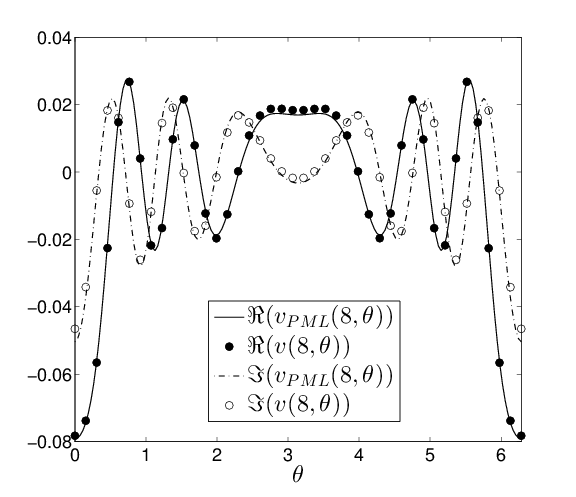}}
\hspace*{-0.5cm}\subfigure[Problem \textbf{PVR 1}: $k=2$]{\includegraphics[width=7.5cm]{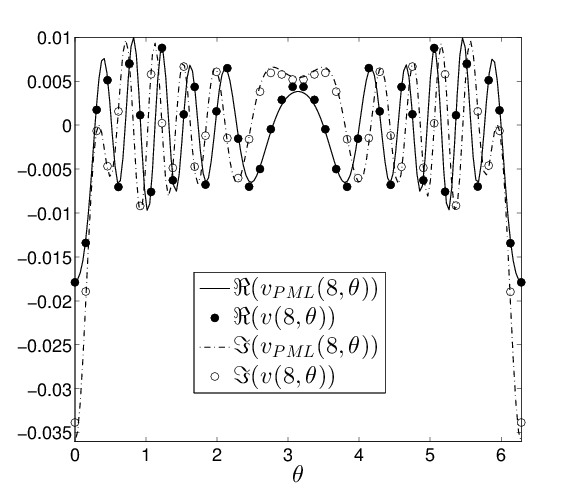}}

\caption{\textbf{Problem PVR 1}: the comparison between the real and imaginary parts of
$v(8,\theta)$ and $v_{PML}(8,\theta)$ for $k=1,2$. The $L^2_{rel}$ error norms are reported in Table 1.}
\label{mono_indvar_k1}
\end{figure}
\end{center}

\begin{center}
\begin{figure}
\hspace*{-1.2cm}\subfigure[Problem \textbf{PVR 2}: $k=1$]{\includegraphics[width=7.5cm]{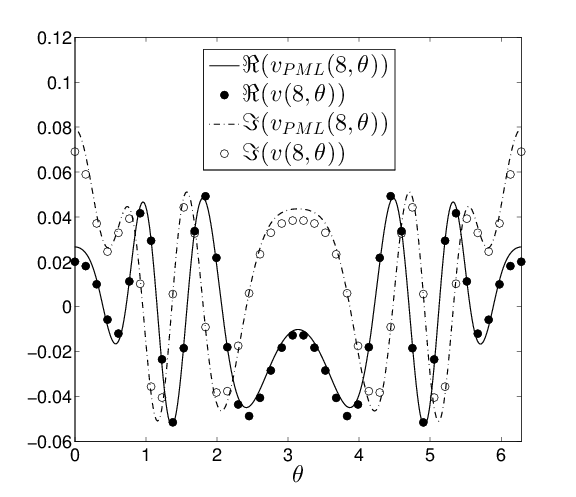}}
\hspace*{-0.5cm}\subfigure[Problem \textbf{PVR 2}: $k=2$]{\includegraphics[width=7.5cm]{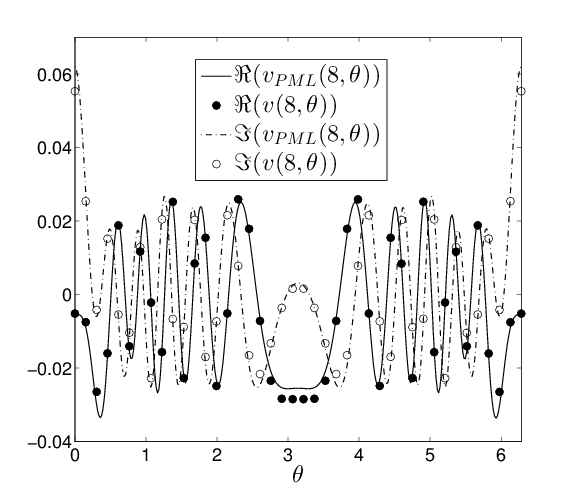}}

\caption{\textbf{Problem PVR 2}: the comparison between the real and imaginary parts of
$v(8,\theta)$ and $v_{PML}(8,\theta)$ for $k=1,2$. The $L^2_{rel}$ error norms are reported in Table 1.}
\label{dipo_indvar_k1}
\end{figure}
\end{center}

The numerical simulations for PVR 1-2 are performed by setting $R=8$, $M_{\rho}=600$ and $M_{\theta}=41$.
For this kind of problems, the exact solution can not be computed but PML can be applied.
Since PML is very efficient and has a very fast convergence, we compare our numerical solution $v(\rho,\theta)$ obtained from the minimization procedure,
to the solution $v_{PML}(\rho,\theta)$ obtained by using a PML method,
where we chose a quadratic turn-on of the PML absorption for a PML layer of thickness
comparable to the wavelength (such choices give usually negligible reflections \cite{OZAJ},\cite{TH}).
In Figures \ref{ppvr1}--\ref{ppvr2}  the real
and imaginary parts
of the numerical solutions $v(\rho,\theta)$ of PVR 1-2 are shown.
In Figures \ref{mono_indvar_k1}-\ref{dipo_indvar_k1}
we show the comparison between the real and imaginary parts of  $v(\rho,\theta)$ and $v_{PML}(\rho,\theta)$ in $\partial B_8$ for $k=1,2$.
In Table 1 we report the $L^2_{rel}$ error norms in $\partial B_8$.  We notice that the two solutions $v$ and $v_{PML}$ are in good agreement
and that the errors in Table 1 are of the same order as the ones in \cite{CGS13}.

\section{Numerical results II: scattering problems} \label{sect 5}
In this Section we present numerical results which are related to some classical scattering problems. In a typical scattering problem one assumes that the total field $u$ can be decomposed into an incident wave $u^i$ and a scattered wave $u^s$, i.e.
\begin{equation*}
u=u^i+u^s.
\end{equation*}
Here, we assume that the index of refraction of the background medium is $n_0(x)$. That is, if we assume that the scatterer is represented by sone domain $D$, then $n(x)=n_0(x)$ for $x\in \RR^2\setminus D$. We consider an incident wave $u^i$ which satisfies
\begin{equation}\label{helm u^i}
\Delta u^i + k^2 n_0(x)^2 u^i = 0, \quad x \in \RR^2.
\end{equation}
Hence, the scattered wave $u^s$ is the solution of
\begin{equation}\label{helm u^s}
\Delta u^s + k^2 n(x)^2 u^s = k^2 [n_0(x)^2-n(x)^2] u^i, \quad x \in \RR^2,
\end{equation}
which satisfies the radiation condition at infinity (like \eqref{rad cond Somm intro} or \eqref{rad con Perth Vega intro}).

In the following, we present some numerical simulations for the case $n_0(x)=1$.

\begin{center}
\begin{figure}
\hspace*{-1.2cm}\subfigure[Problem \textbf{PCR 1}: $k=1$]{\includegraphics[width=7.5cm]{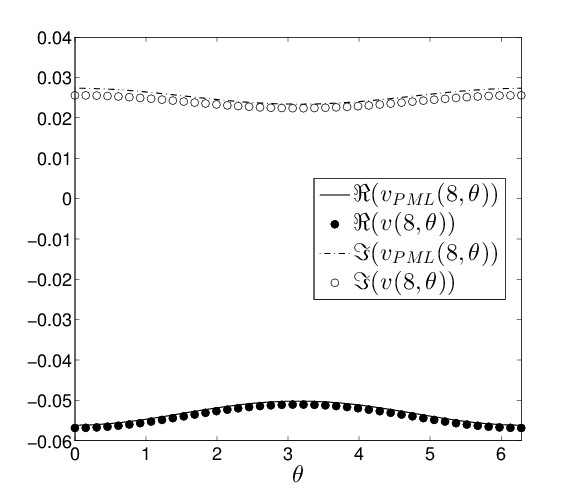}}
\hspace*{-0.5cm}\subfigure[Problem \textbf{PCR 1}: $k=2$]{\includegraphics[width=7.5cm]{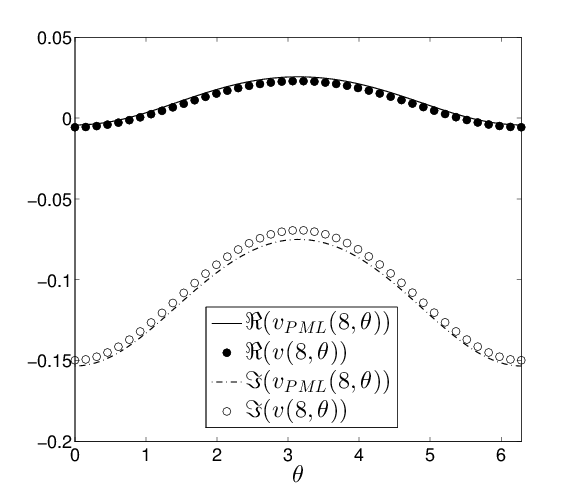}}
\caption{\textbf{Problem PCR 1}: the comparison between the real and imaginary parts of
$v(8,\theta)$ and $v_{PML}(8,\theta)$ for $k=1,2$. The $L^2_{rel}$ error norms are reported in Table 1.}
\label{comp_disco_centrk1}
\end{figure}
\end{center}

\begin{center}
\begin{figure}
\hspace*{-1.2cm}\subfigure[Problem \textbf{PCR 2}: $k=1$]{\includegraphics[width=7.5cm]{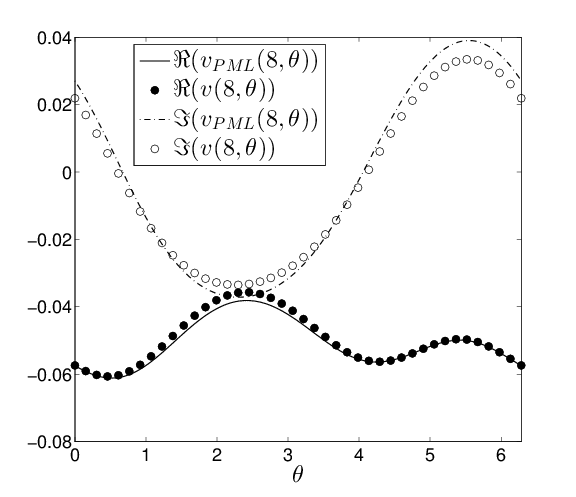}}
\hspace*{-0.5cm}\subfigure[Problem \textbf{PCR 2}: $k=2$]{\includegraphics[width=7.5cm]{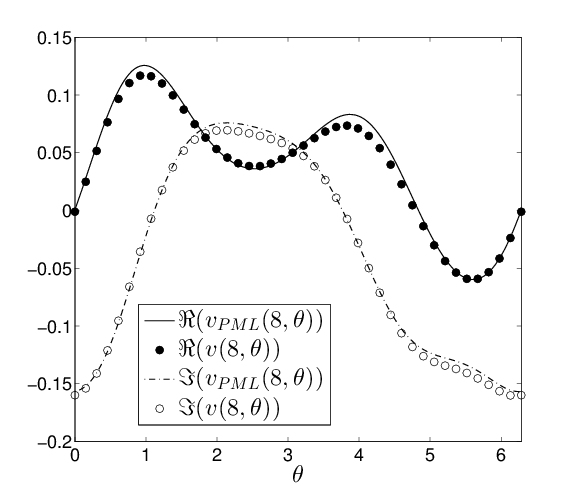}}

\caption{\textbf{Problem PCR 2}: the comparison between the real and imaginary parts of
$v(8,\theta)$ and $v_{PML}(8,\theta)$ for $k=1,2$. The $L^2_{rel}$ error norms are reported in Table 1.}
\label{comp_disco_deck2}
\end{figure}
\end{center}

\begin{center}
\begin{figure}
\vspace*{-0cm}
\hspace*{-1.2cm}\subfigure[Problem \textbf{PCR 4}: $k=1$]{\includegraphics[width=7.5cm]{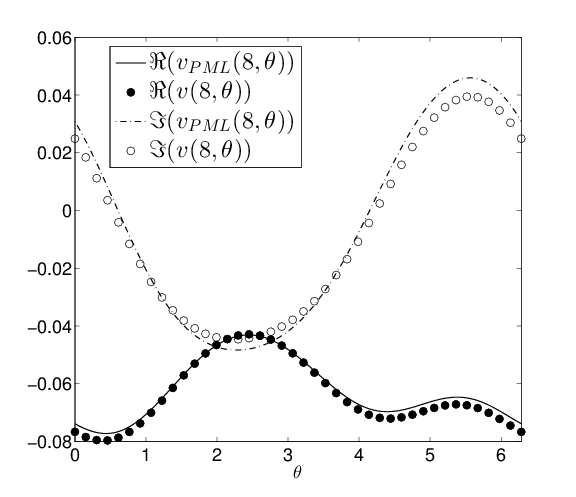}}
\hspace*{-0.5cm}\subfigure[Problem \textbf{PCR 4}: $k=2$]{\includegraphics[width=7.5cm]{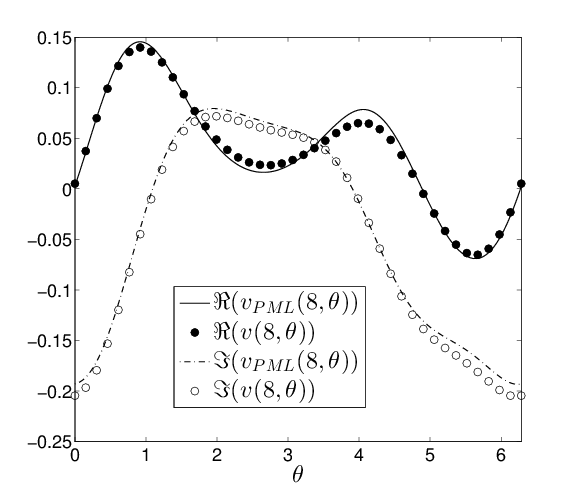}}

\caption{\textbf{Problem PCR 4}: the comparison between the real and imaginary parts of
$v(8,\theta)$ and $v_{PML}(8,\theta)$ for $k=1,2$. The $L^2_{rel}$ error norms are reported in Table 1.}
\label{comp_qua_deck2}
\end{figure}
\end{center}

\paragraph{Background medium with constant index of refraction.} We assume that the background medium has constant index of refraction
$n_0(x)=1$. Let $B_r(p)$ and $Q_r(p)$ be the ball of radius $r$ centered at $x_0$ and the square of side $2r$ centered at $p \in \RR^2$,
respectively, and assume that the incident wave is
$$u^i(x_1,x_2)=\text{e}^{ikx_1}.$$
We set $r=0.5$, $p=(0.5,0.5)$, and consider problem \eqref{helm u^s} by choosing several indexes of refraction:
\begin{description}
 \item[PCR 1]: $n(x)^2= 1+\chi_{B_r(O)}(x)$,
 \item[PCR 2]: $n(x)^2= 1+\chi_{B_r(p)}(x)$,
 \item[PCR 3]: $n(x)^2= 1+\chi_{Q_r(O)}(x)$,
 \item[PCR 4]: $n(x)^2= 1+\chi_{Q_r(p)}(x)$,
\end{description}
where $\chi$ is the characteristic function.
The numerical simulations for PCR 1-4 are performed by setting $R=8$, $M_{\rho}=600$ and $M_{\theta}=41$.
For each problem PCR 1-4, we denote by $v_{PML}(\rho,\theta)$
the numerical solution obtained by using a PML method; also for these problems  we have chosen a quadratic turn-on of the PML absorption for a PML
layer of thickness comparable to the wavelength.
We tested our results by comparing the numerical
solution $v(\rho,\theta)$ , to the solution $v_{PML}(\rho,\theta)$. In Figures \ref{comp_disco_centrk1}--\ref{comp_qua_deck2} we compare the real and imaginary parts of
$v(\rho,\theta)$ and $v_{PML}(\rho,\theta)$ at $\rho=8$ and for $k=1,2$. In Table 1 we report the $L^2_{rel}$ error norms in $\partial B_8$.
Also in these cases the two solutions $v$ and $v_{PML}$ are in good agreement
and  the errors in Table 1 are still of the same order as the ones in \cite{CGS13} and of the problems PVR 1-2.
We observe that the errors are generally greater for the cases in which
the obstacle $D$ is not centered at the origin (PCR 2 e PCR 4).

\begin{center}
\vspace*{-2cm}\begin{figure}
\hspace*{-1.2cm}\subfigure[Problem \textbf{PCR 3}: $k=1$]{\includegraphics[width=7.5cm]{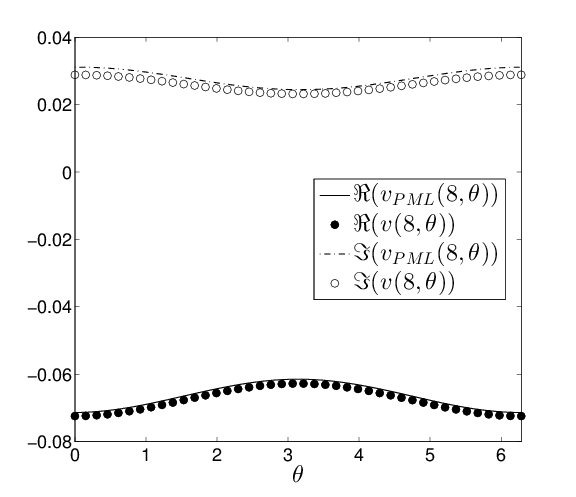}}
\hspace*{-0.5cm}\subfigure[Problem \textbf{PCR 3}: $k=2$]{\includegraphics[width=7.5cm]{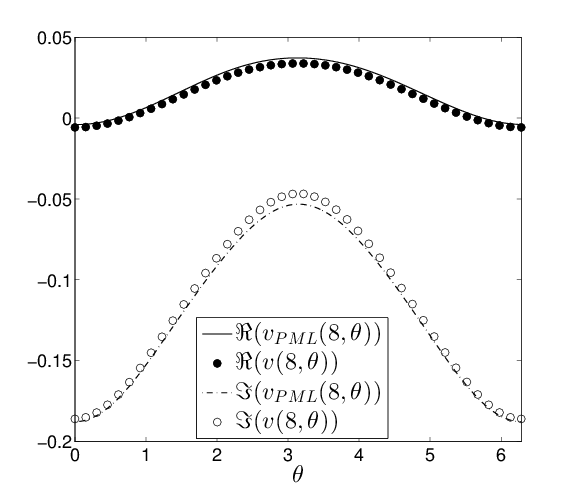}}

\caption{\textbf{Problem PCR 3}: the comparison between the real and imaginary parts of
$v(8,\theta)$ and $v_{PML}(8,\theta)$ for $k=1,2$. The $L^2_{rel}$ error norms are reported in Table 1.}
\label{comp_qua_centrk1}
\end{figure}
\end{center}

\paragraph{Unbounded perturbation.}
Here, we consider a case in which the methods available in literature are of difficult application.
In particular we consider a medium with constant background index of refraction and a perturbation which extends to infinity and has some discontinuity. In particular, we consider the following two scattering problems for \eqref{helm u^s} with $u^i(x)=e^{ikx_1}$:
\begin{description}
 \item[PW 1]: $n(x)^2= 1+\frac{1}{1+|x|}\chi_{E}(x)$,
 \item[PW 2]: $n(x)^2= 1+\frac{1}{1+|(x_1,x_2-2)|}\chi_{E}(x)$,
\end{description}
where $E=\{ (x_1,x_2) : \lvert x_1-\sin(x_2) \rvert \leq 0.25 \}$.

\begin{center}
\begin{figure}
\hspace*{-1.2cm}\subfigure[Problem \textbf{PW 1}, real part: $k=1$]{\includegraphics[width=7.3cm]{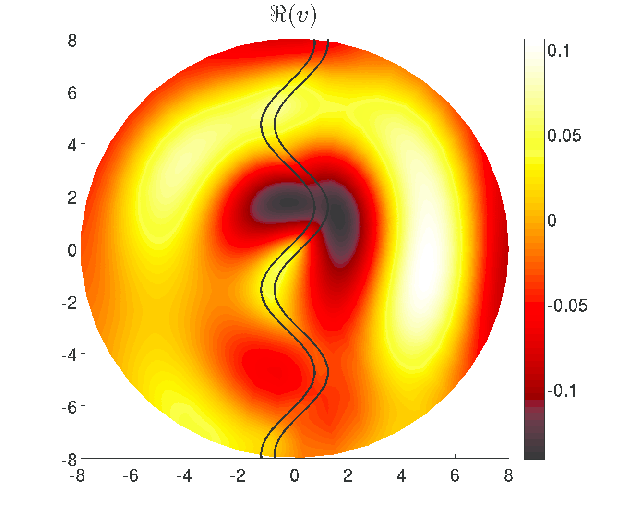}}
\hspace*{-0.5cm}\subfigure[Problem \textbf{PW 1}, imaginary part: $k=1$]{\includegraphics[width=7.3cm]{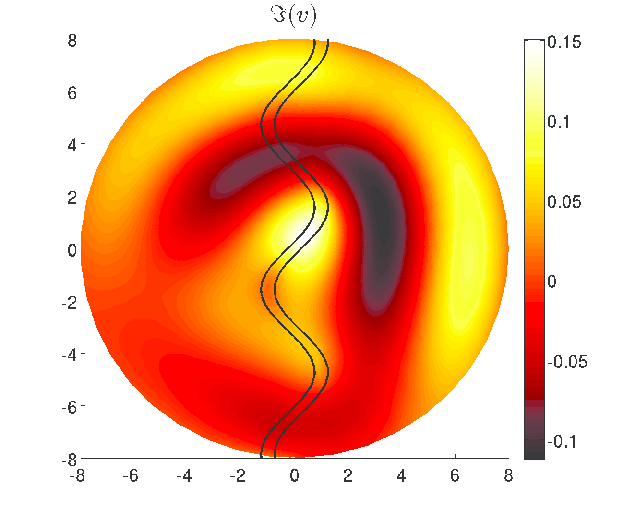}}
\hspace*{-1.2cm}\subfigure[Problem \textbf{PW 1}, real part: $k=2$]{\includegraphics[width=7.3cm]{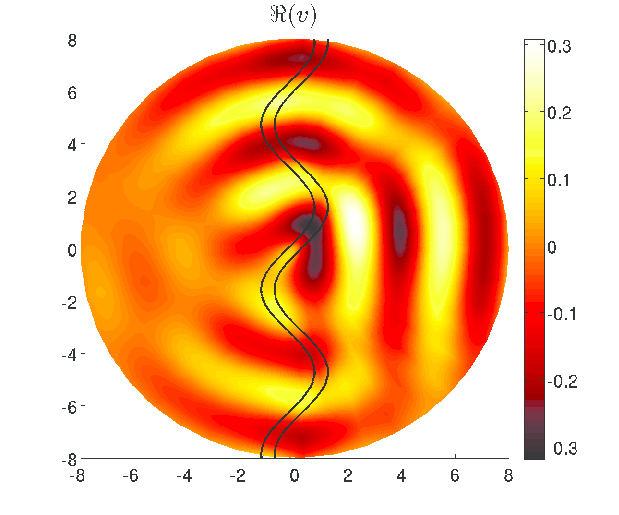}}
\hspace*{-0.5cm}\subfigure[Problem \textbf{PW 1}, imaginary part: $k=2$]{\includegraphics[width=7.3cm]{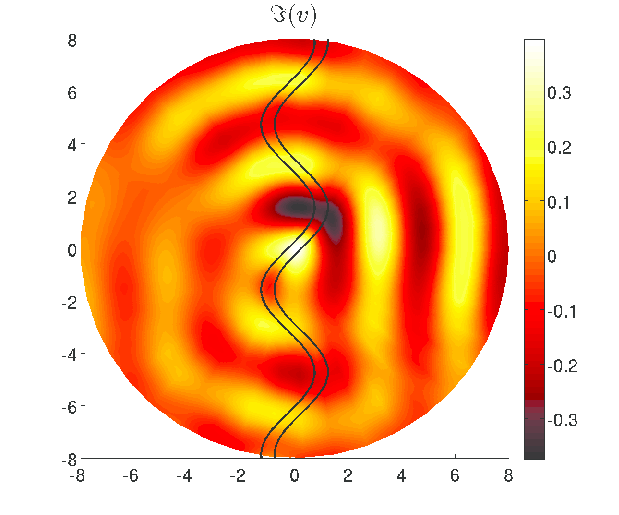}}
\caption{\textbf{Problem PW 1:} real and imaginary part of the numerical solution for $k=1,2$. The black lines delimit the support of the set $E$.}
\label{ppw1}
\end{figure}
\end{center}

\begin{center}
\begin{figure}
\hspace*{-1.2cm}\subfigure[Problem \textbf{PW 2}, real part: $k=1$]{\includegraphics[width=7.3cm]{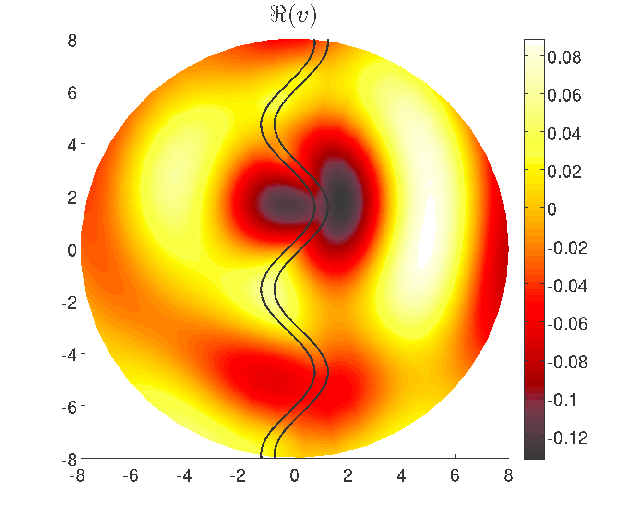}}
\hspace*{-0.5cm}\subfigure[Problem \textbf{PW 2}, imaginary part: $k=1$]{\includegraphics[width=7.3cm]{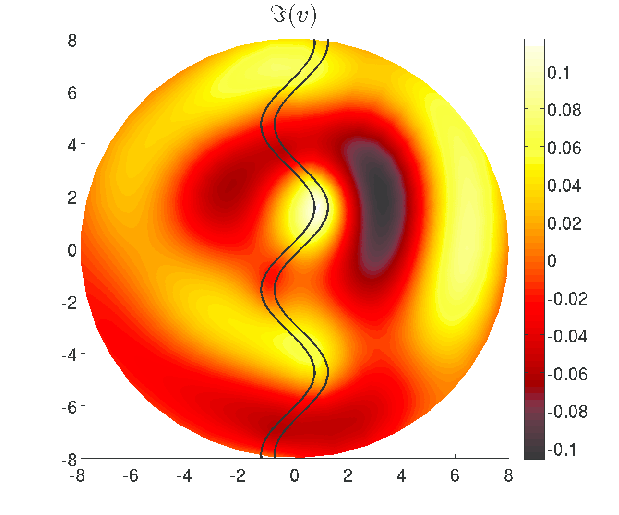}}
\hspace*{-1.2cm}\subfigure[Problem \textbf{PW 2}, real part: $k=2$]{\includegraphics[width=7.3cm]{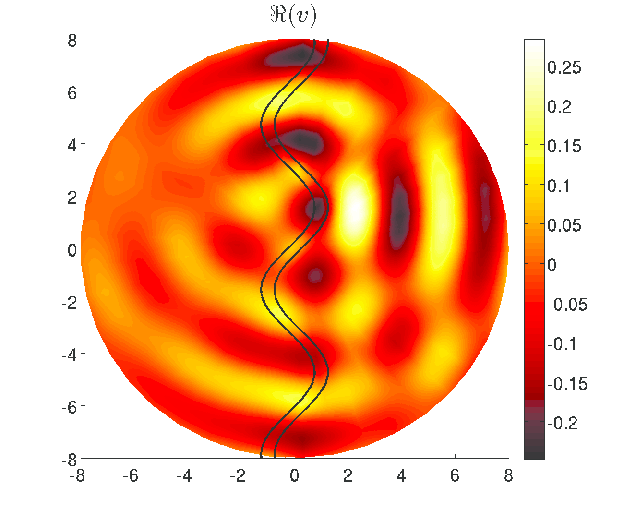}}
\hspace*{-0.5cm}\subfigure[Problem \textbf{PW 2}, imaginary part: $k=2$]{\includegraphics[width=7.3cm]{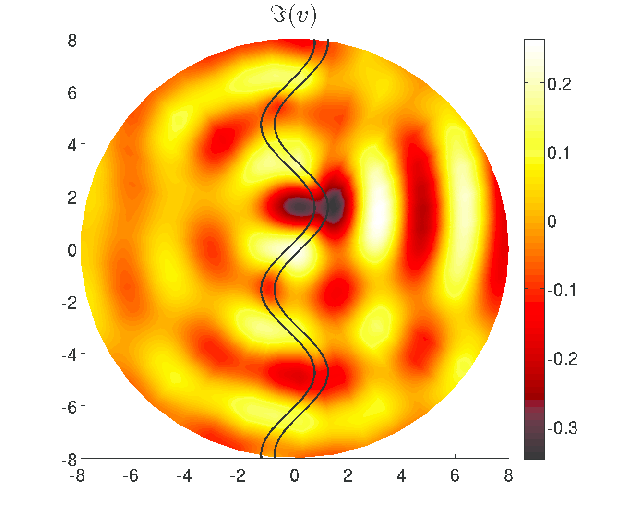}}
\caption{\textbf{Problem PW 1:} real and imaginary part of the numerical solution  for $k=1,2$. The black lines delimit the support of the set $E$.}
\label{ppw2}
\end{figure}
\end{center}

Since the index of refraction is not constant outside any compact region,
we do not know the Dirichlet-to-Neumann map. Moreover,
the index of refraction can not be continued analytically in the direction orthogonal
to the computational domain and so also the PML method seems to be of non-standard application.
As proved in \cite{CGS13}, our approach still works for this kind of problem. We present our numerical results in Figures \ref{ppw1}-\ref{ppw2} where the real
and imaginary parts
of the numerical solutions $v(\rho,\theta)$ of PWR 1-2 are shown.

\section{Conclusions}
We considered a Fourier-Chebyshev collocation method for studying a constrained optimization problem which is related to a new approach to
the problem of transparent boundary conditions for the Helmholtz equation in unbounded domains (see \cite{CGS13}).
We gave numerical evidence of an estimate available in literature \cite{Ci3} which gives the rate of convergence of the numerical
scheme at a fixed frequency $k$. We studied numerically the problem at low and mid-high frequencies and show that the minimization problem improves for $k$ large. However, due to the large computational complexity of the numerical scheme, we were not able to test our approach in the (very) high frequency regime. In the low frequency regime, we observed that the convergence of the numerical results depends on the source function: for some they deteriorate while improve for some others. We believe that the numerical studies presented in this paper give a hint on the convergence properties of the algorithm also in more general settings.

\begin{table}
\begin{center}
    \begin{tabular}{ccc}
    \hline
    $Problem$ & $k$ &$L^2_{rel}$ error  \\
\hline
PCR 1 & 1 & 0.027 \\
PCR 1 & 2 & 0.045 \\
PCR 2 & 1 & 0.073 \\
PCR 2 & 2 & 0.061 \\
PCR 3 & 1 & 0.029 \\
PCR 3 & 2 & 0.042 \\
PCR 4 & 1 & 0.066 \\
PCR 4 & 2 & 0.073 \\
PVR 1 & 1 & 0.049 \\
PVR 1 & 2 & 0.061 \\
PVR 2 & 1 & 0.097 \\
PVR 2 & 2 & 0.099 \\
\hline
\end{tabular}
\end{center}
\caption{Table of $L^2_{rel}$ errors on $\partial B_8$ between the solution $v$ of \eqref{minimiz problem PV sect 2} and the related PML solution $v_{PML}$.}
\label{tableconv2}
\end{table}




%







\end{document}